\documentstyle[12pt]{article}

  \newcommand{\const}{\rm const}
  \newcommand{\Var}{\rm Var}
  
  \newcommand{\Law}{\rm Law}

  \newcommand{\Cov}{\rm Cov}
  \newcommand{\Ent}{\rm  Ent}

\begin{document}

 \begin{center}

 {\bf  Method Monte - Carlo for solving of non - linear integral equations. }

 \vspace{5mm}

 {\bf M.R.Formica, \ E.Ostrovsky, \ L.Sirota. }

 \end{center}

\vspace{4mm}

 \ Universit\`{a} degli Studi di Napoli Parthenope, via Generale Parisi 13, Palazzo Pacanowsky, 80132,
Napoli, Italy. \\

e-mail: mara.formica@uniparthenope.it \\

\vspace{3mm}

 \ Israel,  Bar - Ilan University, department  of Mathematic and Statistics, 59200, \\

\vspace{3mm}

e-mails:\\
eugostrovsky@list.ru \\
sirota3@bezeqint.net \\

\vspace{4mm}

\begin{center}

 {\bf Abstract} \\

 \vspace{4mm}

\end{center}

  \ We offer in this short report a simple Monte - Carlo method for solving a well - posed non - linear integral equations of
  second  Fredholm's and Volterra's type and built a confidence region for solution in an uniform norm, applying the grounded
  Central Limit Theorem in the Banach space of continuous functions. \par

   \ We prove that the rate of convergence our method coincides  with the classical one.\par

 \vspace{4mm}

\begin{center}

 {\bf Key words and phrases.} \\

 \vspace{3mm}

\end{center}

 \  Non - linear  integral equation of the  Fredholm's and Volterra's second type, confidence region, uniform norm, Monte - Carlo method,  metric
 spaces, metric entropy, compactness, Borelian distribution, weak convergence and compactness, H\"older's condition, kernel, distance, Mittag - Leffler's  function,
 contraction principle, partition,  random process  and field (r.f.), Law of Large Numbers  (LLN), Banach space,  Lipschitz  condition,
 convergence, asymptotical approach, Central Limit Theorem (CLT),  Banach space of continuous functions, Gaussian distribution, iterations,  estimations,
 well posedness, approximation, speed of convergence,  Lebesgue - Riesz  spaces, moments, spectral radius,  modulus of the uniform continuity.\\

\vspace{5mm}

\section{Definitions. Statement of problem. Notations. Fredholm's equations.}

\vspace{5mm}

 \ Let $ \ (T =\{t\}, B,\mu, d) \ $ be compact (complete)  metric probabilistic measure space: $ \ \mu(T) = 1, \ $ equipped with a  non - trivial (bounded) distance
 function $ \ d = d(t,s), \ t,s \in T. \ $  It will be presumed as usually that the measure $ \ \mu \ $ is Borelian.\par

 \  We  consider in this preprint  the following non - linear integral equation of the second  Fredholm's type

\begin{equation} \label{equation nonl}
x(t) = f(t) + \int_T K(t,s, x(s)) \ \mu(ds).
\end{equation}

 \ The function $ \ f: T \to R \ $ in (\ref{equation nonl})   will be presumed to be  continuous, $ \  K = K(t,s,z), \ t,s \in T, \ z \in R  \ $ is common continuous kernel,
 non - linear  in general case relative  the last variable $ \ z; \   \ x = x(t), \ x: T \to R \ $ is unknown function. \par

\vspace{3mm}

  \ The case of Volterra's  second type of integral equation will be considered in a penultimate section. \par

\ We agree to take henceforth to write

$$
\int g(s) \ d \mu = \int_T g(s) \ \mu(ds)
$$
for any measurable function $ \ g: T \to R. \ $ \par

\vspace{3mm}

 \ We will deal with the space of all continuous numerical valued functions $ \ C(T,d) = C(T) \ $   equipped with an ordinary uniform norm

$$
||x|| = ||x||C(T) \stackrel{def}{=}  \sup_{t \in T} |x(t)| = \max_{t \in T}|x(t)|.
$$

 \ The case of the Lebesgue - Riesz spaces $ \ L_p(T,\mu) \ $ is considered in many articles,  in particular in  \cite{Karoui}.

\vspace{4mm}

 \ We must introduce in this regards some additional notations and assumptions.  Namely, let $ \ f = f(t) \ $ be (bounded) continuous function:
 $ \ f(\cdot) \in C(T). \ $  Further, impose the following important condition on the kernel function $ \ K = K(t,s,z): \ $ it will be bounded,  common continuous
and  satisfies relative the third argument  Lipschitz  condition  with constant less than 1:

\begin{equation} \label{rho condition}
\exists \rho = \const \in (0,1) \ \Rightarrow |K(t,s,z_1) - K(t,s,z_2)| \le \rho |z_1 - z_2|.
\end{equation}

 \ It is known, see e.g. \cite{Kantorowisz Akilov}, chapter 16, section 4;  \cite{Krasnoselskii} that the source problem   (\ref{equation nonl}) is well - posed:
 the {\it continuous} solution $ \ x = x(t) \ $ there exists, is unique and continuously depended on the datum $ \ f(\cdot), \ K(\cdot,\cdot,\cdot),  \ $ of course,
 subject to the condition (\ref{rho condition}). See also the classical monographs \cite{Dunford Schwartz1}, chapters 2,3;   \cite{Dunford Schwartz2}, chapter 1. \par

 \ Moreover,  the solution $ \ x = x(t) \ $ may be obtained  in particular as an {\it uniform} limit as $ \ n \to \infty \ $ of the following recursion:
   $ \ x_0(t) := f(t), \ $

 \begin{equation} \label{recursion}
x_{n+1}(t) := f(t) + \int_T K(t,s,x_n(s)) \ \mu(ds), \ n = 0.1,2,\ldots.
 \end{equation}
 \ Briefly:

\begin{equation} \label{recursion briefly}
x_{n+1}(t) :=  R[f,K,x_n](t); \
\end{equation}
where

\begin{equation} \label{def recursion}
 R[f,K,x_n](t) =  \int_T K(t,s,x_n(s)) \ \mu(ds), \ n = 0.1,2,\ldots.
\end{equation}

\vspace{3mm}

\ An error estimate:

\begin{equation} \label{err est}
||x - x_m|| \le ||x_1 - x_0|| \cdot \frac{\rho^m}{1 - \rho}. \ m = 1,2,\ldots,
\end{equation}
contraction principle. \par

 \ We will agree to write for certain random field (r.f.)  $ \ \Delta(t) \ $

$$
\Law(\Delta(\cdot)) = N(a(t), \ R(t,s)), \ t,s \in T
$$
if the r.f.  $ \ \Delta(t) \ $ has a Gaussian (normal) distribution with parameters

$$
{\bf E}\Delta(t) = a(t), \ \Cov(\Delta(t), \Delta(s)) = R(t,s).
$$
 \  As a rule in this report $ \ a(t) = 0; \ - \ $  the so - called centered case. \par

  \ Of course, the function  $ \ R = R(t,s) \ $ must be finite, symmetrical and non - negative definite. \par
 \ Further, let $ \ \{\Delta_N \} = \{\Delta_N(t)\}, \ t \in T, \ N \ = 1,2,\ldots \ $ be a {\it sequence}  of continuous a.e. random fields. \par

 \ We will agree to write for the sequence of  random fields $ \ \Delta_N(t) \ $ converges weakly (in distribution) in this space  $ \ C(T,d) = C(T), \ $
 as $ \ N \to \infty, \ $ to some r.f. $ \ \Delta(t), \ $

$$
\Law(\Delta_N(\cdot)) \stackrel{dist}{\to}  \Law(\Delta(\cdot)).
$$
  \ As a rule, the limiting r.f. $ \ \Delta \ $ will be Gaussian and centered. \par
 \  We  recall further some facts about this convergence (Prokhorov - Skorokhod theory). \par

\vspace{5mm}

\section{ Description of method. Convergence.}

\vspace{5mm}

 \   But  the question arises: how to calculate the integrals holding in  (\ref{recursion}) integrals? \par

\vspace{4mm}

 \ {\bf Our claim in this report is  to offer the consistent Monte - Carlo method, more precisely, a depending trials method,  for these integral computations,
and to built an asymptotical as}  $ \ N \to \infty \ $  {\bf confidence domain  for solution in an uniform norm. } \par

\vspace{4mm}

 \ We will prove that the rate of  convergence as $ \ N \to \infty \ $  of this method  in the uniform norm for the $ \ x_m(t) \ $ calculation is equal to
 $ \  1/\sqrt{N}, \ $ where $ \ N \ $ is  an amount of common elapsed random numbers having  the distribution $ \ \mu. \ $ \par
  \ More precisely, we will ground that the classical normed $ \  \sqrt{N}, \ $ where $ \ N \ $ denotes the amount of all elapsed random variables,
  deviation between Monte - Carlo approximation and deterministic one satisfies the Central Limit Theorem
 (CLT) in the space of all continuous functions. \par

 \vspace{4mm}

 \ In detail. Let $ \ \vec{\xi} =   \vec{\xi}[N] = \{ \xi(1),\xi(2),\xi(3),\ldots,\xi(N)  \}, \ \xi = \xi(1) \ $ be a $ \ N - \ $
 tuple of independent random variables (r.v.-s.)
having the distribution $ \ \mu: \ $

$$
{\bf P} (\xi(i) \in A) = \mu(A)
$$
and $ \ N \ $ is some "great" integer number $ \ N >> 1. \ $ \par

 \ Introduce a set of integer numbers

\begin{equation} \label{set N}
S = S(N) = \{ \ 1,2,3,\ldots, N-1,N \ \},
\end{equation}
and introduce also some  its partition $  \ Q = \{Q_k\}, \ k = 1,2,  \ldots,m \ $ as follows

\begin{equation} \label{partition one}
Q(1) = \{1,2,\ldots,n(1)\}; \ Q(2) = \{n(1)+1,n(1) + 2, \ldots, n(2)\}, \ \ldots,
\end{equation}

\begin{equation} \label{partition general}
 Q(m-1) = \{n(m-1)+1, n(m-1)+2. \ldots, n(m)  \}.
\end{equation}

\ Of course, $ \ \{n(j)\} \ $ are positive integer numbers and such that
$ \  1 < n(1) < n(2) \ldots < n(m-1) < n(m) = N. \ $ \par

\vspace{4mm}

 \ Introduce the following integer vector depending  (in general case) on the value $ \ N, \ $
 (and of course on the number of iterations $ \ m) \ $  of differences
 $ \ \vec{q} =  \vec{q[N]} = \{ q(k)  \}, \ k = 1,2,\ldots,m: \ $

\begin{equation} \label{vec q}
q(1)= q[N](1), \hspace{4mm} q(2) = q[N](2) = n(2) - n(1),
\end{equation}

\begin{equation} \label{vec qm}
q(3) = q[N](3) = n(3) - n(2), \hspace{4mm}  \ldots,   q(m) = q[N](m) = n(m) - n(m-1).
\end{equation}

 \ Evidently, the numbers $ \ q[N](k) \ $ are positive integer and

\begin{equation} \label{q norming}
\sum_{k=1}^m q[N](k) = n(m) = N.
\end{equation}

\vspace{4mm}

 \ {\bf We impose in the sequel throughout all the report the following important conditions on the introduced vector either  }

\begin{equation} \label{comdit q m -1}
\lim_{N \to \infty} q[N](k) = \infty, \  \underline{\lim}_{N \to \infty} q[N](m)/N   > 0;
\end{equation}

\vspace{3mm}

 or sometimes a more strong condition

\begin{equation} \label{ comdition q m }
\lim_{N \to \infty} q[N](k) = \infty, \  \lim_{N \to \infty} q[N](m)/N  = 1.
\end{equation}

\vspace{5mm}

 \ Let us introduce also a following vector  depending on the variable $ \ N  \ $ and on the value $ \ m: \ $

\begin{equation} \label{vec gamma}
  \gamma = \vec{\gamma} = \vec{\gamma}_N = \{ \ \gamma_N(i) \ \}, \ i = 1,2,\ldots,m, \ \dim \gamma = m,
\end{equation}

 \begin{equation}  \label{n m vector}
\gamma(1) = \gamma_N(1) = n(1)/N, \ \gamma(2)= \gamma_N(2) = [n(2) - n(1)] /N,
\end{equation}

\begin{equation} \label{lastnm}
\gamma(3) = \gamma_N(m) = [n(3) -n(2)]/N, \ldots, \   \gamma(m) =  \gamma_N(m) = [n(m) - n(m-1)]/N.
 \end{equation}
 \ Evidently, $ \ \gamma_N(k) \in (0,1) \ $ and

 $$
  \sum_{k=1}^m \gamma(k) =  \sum_{k=1}^m \gamma_N(k) = 1.
 $$

\vspace{4mm}

 \  {\bf  We introduce  also the condition that as}  $ \ N \to \infty \  \Rightarrow   $

\begin{equation} \label{condition gamma}
0 < \underline{\lim}_{N \to \infty} \min_k \gamma_N(k) \le \overline{\lim}_{N \to \infty} \max_k \gamma_N(k) < 1.
\end{equation}

\vspace{4mm}

\  In particular, the sequence $ \ \gamma_N(k) \ $ can be selected as

\begin{equation} \label{select gamma}
\gamma_N(k) = 1/m, \ k = 1,2,\ldots,m
\end{equation}

 \ Evidently, the condition (\ref{condition gamma}) entails one  (\ref{comdit q m -1}). \par

\vspace{5mm}

 \ The {\it optimal choice} of the variables $ \ q[N](k), \ k = 1,2,\ldots, m \ $ will be clarified below. \par

\vspace{5mm}

 \ Let us offer the following {\it iterative} Monte - Carlo procedure   $ \ x^{0}_0(t)  = x_0(t) := f(t), \ $
and for the values $ \ k  = 1,2, 3,4,\ldots, m-1,m, \ $

 \begin{equation} \label{first rand recursion}
x_1^1(t) = x_1^{1(N)}(t) = f(t) +  \frac{1}{n(1)} \sum_{i=1}^{n(1)} K(t,\xi(i),x_0^0(\xi(i))),
 \end{equation}

 \begin{equation} \label{ second rand recursion}
x_2^2(t) = x_2^{2(N)}(t) = f(t) +  \frac{1}{n(2) - n(1)} \sum_{i=n(1)+1}^{n(2)} K(t,\xi(i),x_1^1(\xi(i))),
 \end{equation}

 \begin{equation} \label{ k  rand recursion}
x_k^k(t) = x_k^{k(N)}(t) = f(t) +  \frac{1}{n(k) - n(k-1)} \sum_{i=n(k-1)+1}^{n(k)} K(t,\xi(i),x_{k-1}^{k-1}(\xi(i))),
 \end{equation}

 \begin{equation} \label{m rand recursion}
x_m^m(t)  = x_m^{m(N)}(t) = f(t) +  \frac{1}{n(m) - n(m-1)} \sum_{i=n(m-1)+1}^{n(m)} K(t,\xi(i),x_{m-1}^{m-1}(\xi(i))),
 \end{equation}
and recall that $ \ n(m) = N. \ $

\vspace{5mm}

 \ {\bf Theorem 2.1.} Suppose in addition that the metric space $ \ (T,d) \ $ is (complete) compact. Then the sequence of random fields
 $ \  x_m^{m(N)}(t)  \ $  converges as $ \ N \to \infty \ $ uniformly with probability one to the approximation $ \ x_m = x_m(t), \ t \in T: \ $

\begin{equation} \label{uniform convergebce}
{\bf P} \left( \ \lim_{N \to \infty} \max_{t \in T}|x_m^m(t) - x_m(t)| \to 0 \  \right) = 1,
\end{equation}
 as well as this proposition holds true for all the previous values $ \ k. \ $ \par

\vspace{4mm}

 \ {\bf Proof.} Note that the Banach space $ \ C(T,d) = C(T) \ $ of all numerical values continuous functions is complete and separable. Therefore one can
  apply the famous Law of Large Numbers  (LLN) in this Banach space, see e.g.  \cite{Fortet Mourier1}, \cite{Fortet Mourier2}, \cite{Ledoux Talagrand}, chapters 4,5.
As long as the kernel $ \  K = K(t,s,z) \ $ is continuous and bounded, it follows from (\ref{first rand recursion}) that

$$
{\bf P} (||x_1^1 - x_1||C(T) \to 0) = 1.
$$
 \ By means of an induction

$$
{\bf P} (||x_k^k - x_k||C(T) \to 0) = 1, \ k = 2,3,\ldots,m.
$$
 \ The proposition of theorem 2.1 follows for the value $ \ k = m. \ $ \par

\vspace{4mm}

\ {\bf Example 2.1.} L:et us consider the problem of computations (multiple, in general case) parametric  integrals  of the form

\begin{equation} \label{param int}
I(t) = \int_T g(t,s) \ \mu(ds), \ t \in T
\end{equation}
by means of the Monte - Carlo method

\begin{equation} \label{dep trials}
I_N(t) = N^{-1} \sum_{i=1}^N g(t, \xi(i)),
\end{equation}
in our notations. \par
 \ This method  appear at first by A.S.Frolov  and N.N.Tchentzov in \cite{Frolov  Tchentzov}, 1962, and was named as
 "depending trials method", as long as the values $ \ I_n(t_1) \ $ and $ \ I_n(t_2), \ t_{1,2} \in T \ $ are in general case correlated. \par
 \ The modification of this method for discontinuous functions is offer in \cite{Grigorjeva  Ostrovsky}, for non - linear Partial Differential Equations
 of Navier - Stokes type - in  \cite{Ostrovsky3}.\par

 \ Assume now that for almost all the values $ \ s \ $ the function $ \ g(t,s) \ $ is continuous relative the  value $ \ t \ $ and is bounded:

$$
\sup_{i,s \in T} |g(t,s)| < \infty.
$$
 \ Define the probability of the uniform convergence

\begin{equation} \label{P probab one}
{\bf P}_{uc}[g] := {\bf P}(\sup_{t \in T}|I_N(t) - I(t)| \to 0).
\end{equation}

  \ Then

\begin{equation} \label{Frolov Tchentzov}
{\bf P}_{uc}[g] = 1.
\end{equation}

\vspace{5mm}

\section{Main result. Investigation of convergence. Asymptotic confidence region.}

\vspace{5mm}

 \ {\sc  Central Limit Theorem in the space of continuous functions. } \par

\vspace{4mm}

 \ For the concrete error estimate, on the other words, for the building of the confidence  domain in the uniform norm
 we need to evaluate the following tail probability (inside the framework of example 2.1)

\begin{equation} \label{key tail probab}
P_N(u) \stackrel{def}{=} {\bf P}(\sqrt{N} \ ||I_N  - I|| > u), \ u = \const \ge 1.
\end{equation}

 \ In order to evaluate the tail probability $ \ P_N(u), \ $ we need to apply the so - called Central Limit Theorem (CLT)
 in the Banach space of all continuous functions $ \ C(T). \ $ \par
  \ Let us recall some used facts from this theory. Lei $ \ \eta_i = \eta_i(t), \ t \in T, \ \eta(t) := \eta_1(t)  \ $ be a sequence of centered
  (mean zero): $ \ {\bf E}\eta_i(t) = 0 \ $ independent identical distributed random fields, having uniformly bounded second moment:

\begin{equation} \label{second moment}
\sigma^2[\eta] := \sup_{t \in T} {\bf E} \eta^2(t) < \infty.
\end{equation}

   Denote the normed sum

\begin{equation} \label{Sn}
S_N(t):= N^{-1/2} \sum_{i=1}^N \eta_i(t), \ N = 1,2,\ldots.
\end{equation}

 \ Evidently, the finite - dimensional distributions of $ \ S_N(t) \ $ converge as $ \ N \to \infty \ $ to ones for Gaussian  centered r.f. $ \ S(t) = S_{\infty}(t) \ $
having at the same covariation  function  $ \ R_{S}(t,s) \ $ as $ \ \eta(t): $

$$
 R_{S}(t,s) = R_{\eta}(t,s) =   R_{S_N}(t,s) =  {\bf E} S(t) S(s) = {\bf E}\eta(t) \ \eta(s), \ t,s \in T.
$$

 \ By definition, the sequence of r.f. $  \ \{\eta_i(t) \}, \ $   or simple  the  {\it individual} r.f. $ \ \eta(t), \ $ satisfies the CLT
 in the (Banach) space $ \ C(T)\ $ iff the sequence of distributions  of $ \ S_N(\cdot) \ $ in the space $ \ C(T) \ $ converges weakly to one for $ \ S(\cdot): \ $  that is, for
 arbitrary bounded continuous functional $ \ G: C(T) \to R \ $

\begin{equation} \label{weak conver}
\lim_{N \to \infty} {\bf E} G(S_N) = {\bf E}G(S).
\end{equation}

 \ If (\ref{weak conver}) there holds, then take place the convergence of correspondent tail functions

 \begin{equation} \label{tail converg}
 \lim_{N \to \infty} {\bf P}(||S_N|| > u) = {\bf P}(||S|| > u), \ u > 0.
 \end{equation}

\vspace{3mm}

\ The asymptotical  behavior as $ \ u \to \infty, \  $ as well as non - asymptotical estimations for  the tail probability
 $ P_{\infty}(u) = P(u), \ $  is investigated in many works, see e.g. \cite{Dmitrovsky}, \cite{Fernique},  \cite{Ostrovsky0}, chapter 3, sections 3.1 - 3.5;
\cite{Piterbarg} etc. Roughly speaking,

$$
\ln P(u) \sim  - \frac{u^2}{2 \max_{t \in T} R(t,t)} =   - \frac{u^2}{2 \ \sigma^2[\eta] }.
$$

\vspace{4mm}

 \ The Cental Limit Theorem (CLT) in the space of continuous functions $ \ C(T,d) \ $ and its applications are devoted many works, e.g.
\cite{Buldygin1}, \cite{Buldygin2}, \cite{Dudley}, \cite{Frolov  Tchentzov}, \cite{Grigorjeva  Ostrovsky}, \cite{Kozachenko  Ostrovsky2},
\cite{Ostrovsky0}, chapter 4, section 4.4; \cite{Ostrovsky2}  etc.\par

  \ Let us quote as an example the following result belonging to G.Pizier  \cite{Pizier}. Define for some value $ \ p \ge 2 \ $ the following  bounded
 semi - distance natural  function

\begin{equation} \label{natur p  dist}
d_p[\eta](t_1, t_2) \stackrel{def}{=} || \eta(t_1) - \eta(t_2)||_p, \ t_1, t_2 \in T.
\end{equation}

\vspace{3mm}

 \ Hereafter the notation $ \ ||\theta||_p \ $ denotes the usually Lebesgue - Riesz $ \ L(p) \ $ norm of the r.v. $ \ \theta: \ $

$$
||\theta||_p = \left[\ {\bf E} |\theta|^p \ \right]^{1/p}.
$$

\vspace{4mm}

\ Denote by $ \ N(T,d,\epsilon) \ $ the minimal number of $ \ d \ $  closed balls covering the whole space $ \ T. \ $ Evidently,
$ \  \forall \epsilon > 0 \ \Rightarrow  N(T,d,\epsilon) < \infty \ $ iff the metric set $ \ (T,d) \ $ is pre - compact set. \par

 \ The value $ \  H(T,d,\epsilon) = \ln  N(T,d,\epsilon) \ $ is named  as metric entropy of the set $ \ T \ $ relative the distance $ \ d. \ $
Some examples of entropy evaluation may be found, e.g. in \cite{Ostrovsky0}, chapter 3, sections 3.1 - 3.3.  In particular, if the set $ T \ $
is closed bounded subset of the Euclidean  space $ \ R^d \ $ equipped with  ordinary Euclidean distance  $ \ |t_1 - t_2| \ $  and for which

\begin{equation} \label{metric estimate}
r(t,s) \asymp |t-s|^{\alpha}. \ \alpha = \const \in (0,1],
\end{equation}
then

\begin{equation} \label{entropy estim}
N(T,r.\epsilon) \asymp \epsilon^{-d/\alpha}, \ \epsilon \in (0,1).
\end{equation}

 \vspace{3mm}

 \ It is known, see the famous work of G.Pizier \cite{Pizier}, that if

\begin{equation} \label{Pizier cond}
\int_0^1 N^{1/p}(T,d_p[\eta], \epsilon) \ d \epsilon < \infty,
\end{equation}
then the sequence of random fields $ \  \{ \eta_i(t) \ \}  \ $ satisfies the CLT in the space $ \ C(T,d_p). \ $ See also a more general proposition in
\cite{Ostrovsky0}, chapter 3, section 3.17. \par

 \ The last condition (\ref{Pizier cond}) is satisfied if for example $ \ T \ $ is closed bounded subset of whole Euclidean space $ \ R^d \ $
with correspondent norm $ \ |t| \ $  and if  $ \exists \beta > 0, \ \exists  C < \infty \ \Rightarrow \ $

\begin{equation} \label{Koplmogorov cond}
{\bf E} |\eta(t) - \eta(s)|^p \le C |t-s|^{d + \beta},
\end{equation}
Kolmogorov - Slutsky condition.  \par

 \ Another version of the CLT in the space of continuous functions may be found in particular in \cite{Jain Marcus}. \par

\vspace{5mm}

 \ Let us return to the source tail function (\ref{key tail probab}). We need to apply mentioned before CLT in the space $ \ C(T). \ $
 Denote

$$
\zeta_N(t) = \sqrt{N} \ \left( I_N(t) - I(t)  \right).
$$

 \ Suppose that the r.f. $ \  \zeta_i(t) = g(t,\xi(i)) - I(t), \ \zeta(t) = \zeta_1(t) = g(t,\xi) - I(t), \ \xi = \xi_1 \ $ satisfies the
 Central Limit Theorem (CLT) in the space $ \ C(T). \ $ This implies by definition that the sequence of
distributions of r.f. $ \  \zeta_N(\cdot) \ $   in the space $ \ C(T) \ $ converges  as $  N \to \infty $   weakly, i.e.
in distribution,  to the  centered continuous Gaussian distributed random field $ \ \zeta(t) = \zeta_{\infty}(t) \ $
having at the same covariation function as  $ \ \zeta_1(t): \  R[\zeta](t_1,t_2) := $

$$
 {\bf E} (\zeta(t_1) - I(t_1))(\zeta(t_2) - I(t_2)) = \int_T g(t_1,s) g(t_2,s) \ \mu(ds) - I(t_1) I(t_2).
$$
 \ Then

\begin{equation} \label{CLT CT}
\lim_{N \to \infty} P_N(u) = P_{\infty}(u) = P(u), \ u > 0,
\end{equation}
where

\begin{equation} \label{Gauss tail}
 P_{\infty}(u) = {\bf P}(\sup_{t \in T} |\zeta_{\infty}(t)| > u).
\end{equation}

\ As a consequence: as $ \ u \to \infty \ $

$$
\ln P(u) \sim  - \frac{u^2}{2 \max_{t \in T} R[\zeta](t,t)}.
$$

\vspace{4mm}

 \ In order to  make sure the CLT for the random field $ \ g(t,\xi) - I(t), \ $ it is sufficient to refer the Pizier
 conditions. Namely, suppose

$$
\beta = \beta^2[g] = \sup_{t \in T} \Var [ \ g(t,\xi) \ ] < \infty,
$$

$$
\exists p \ge 2 \ \Rightarrow \ \sup_{t \in T} ||g(t,\xi)||_p < \infty,
$$

and that

\begin{equation} \label{Pizier two}
\int_0^1 N^{1/p} [T, d_p[g], \epsilon] d \epsilon < \infty,
\end{equation}

where the distance function $ \  d_p[g](t_1,t_2) \ $  is defined above

$$
d_p[g](t_1,t_2) := || \ [g(t_1,\xi) - g(t_2,\xi)] - [ I(t_1) - I(t_2)]  \  ||_p,
$$

 \ As was notes above, the condition (\ref{Pizier two}) is satisfied if $ \ T \ $ is bounded closed subset of whole space $ \ R^d \ $
 with ordinary Euclidean norm  $ \ |t|, \ $ and if

\begin{equation} \label{Koplmogorov Slutsky cond}
{\bf E} |\ g(t,\xi) - g(s,\xi)|^p \le C |t-s|^{d + \theta}, \exists \theta > 0.
\end{equation}

\vspace{4mm}

 \ It is clear that if the metric $ \ d_p[g] (\cdot,\cdot) \ $ is continuous relative the source one $ \ d(\cdot,\cdot), \ $ then the CLT
in the space $ \ C(T,d[g]_p)  \ $ entails one in $ \ C(T,d). \ $ \par

\vspace{5mm}

  \ {\sc   Preliminary considerations.} \par

\vspace{4mm}

 \ Let us investigate the  first random approximation (\ref{first rand recursion}). Assume that the r.f $ \  K(t,\xi_i, f(\xi_i)) \ $
 satisfies the CLT  in the space $ \ C(T,d); \ $ then one can write (approximately)

$$
x_1^1(t) =  f(t) +  \frac{1}{n(1)} \sum_{i= 1}^{n(1)} K(t,\xi_i, x_0^0(\xi_i)) =
$$

$$
 f(t) +  \frac{1}{n(1)} \sum_{i= 1}^{n(1)} K(t,\xi_i,  f(\xi_i)) =
$$

 \begin{equation} \label{first rand recursion 1}
 x_1(t) +  \frac{1}{\sqrt{ n(1)}} \zeta_1(t) = x_1(t) + \frac{1}{\sqrt{ q[N](1)}} \zeta_1(t),
 \end{equation}
where

\begin{equation} \label{first 1asympt1}
\Law(\zeta_1 (\cdot)) \stackrel{dist}{\to} N \left(0, \ R_1[K,f](t_1.t_2) \ \right).
\end{equation}
and

$$
R_1[K,f](t_1,t_2) = R_1(t,s)  :=  \int_T \ K(t_1,s,f(s)) \ K(t_2,s,f(s)) \mu(ds) -
$$

$$
 \int_T K(t_1,s,f(s))\mu(ds) \ \int_T K(t_2,s,f(s)) \ \mu(ds).
$$

 \ Let us suppose temporarily

\begin{equation} \label{first  normal asympt1}
\Law(\zeta_1(\cdot)) = N \left(0, R_1[f](t_1,t_2) \ \right).
\end{equation}

\vspace{4mm}

 \ The case of the second iteration is more complicated. Let us impose in addition  the following
condition  on the kernel function $ \ K = K(t,s,z): $

$$
\ \exists \lambda = \lambda(t,s), \ \theta = \theta(t,s; z,v), \ \delta = \const  \in (0,1]
$$
such that
$$
\sup_{t,s \in T} |\lambda(t,s)| < \infty, \ \sup_{t,s \in T, \ z,v \in R } |\theta(t,s; z,v)| < \infty, \ \Rightarrow
$$

\begin{equation} \label{smooth  kernel cond}
K(t,s,z) - K(t,s,v) = \lambda(t,s) \ (z-v) + \theta(t,s;z,v) |z-v|^{1 + \delta}.
\end{equation}

 \ We derive under this condition (\ref{smooth  kernel cond}) as $ \ N \to \infty \ $

$$
x_2^2(t) = f(t) + \frac{1}{n(2) - n(1)} \sum_{i=n(1) + 1}^{n(2)} K(t,\xi_i, x_1^1(\xi_i))  \sim  f(t) +
$$

$$
\frac{1}{q[N](1) } \sum _{i=n(1) + 1}^{n(2)} K\left[(t,\xi_i, x_1(\xi_i) + [n(1)]^{-1/2} \ \zeta_1(\xi_i) )\ \right] \sim
$$

$$
f(t) + \int_T K(t,s, x_1(s)) \ \mu(ds) + \frac{1}{\sqrt{n(2) - n(1)}} \zeta_2(t) =
$$

$$
x_2(t) + (q[N](2) N)^{-1/2}\zeta_2(t),
$$
where $ \ \zeta_2(t) = \zeta_2[K,f;N](t) \ $ is (approximately) Gaussian centered continuous random field with covariation function

$$
 R_2(t_1,t_2) = R_2[f,K](t_1,t_2) = \int_T K(t_1,s,x_1(s)) \  K(t_2,s,x_1(s)) \mu(ds) -
$$

$$
\int_T K(t_1,s,x_1(s)) \ \mu(ds) \cdot \int_T K(t_2,s,x_1(s)) \ \mu(ds).
$$

\vspace{4mm}

 \ We find quite  analogously for the values $ \ k = 3,4,\ldots,m, \ $ especially for the "final" value $ \ k = m \ $
 the approximate representations

\begin{equation} \label{final representation}
x_k^k(t) = x_k(t) + (q[N](k))^{-1/2} \zeta_k(t),
\end{equation}
where $ \ \zeta_2(t) = \zeta_2[K,f;N](t) \ $ is (approximately) Gaussian centered continuous random field with covariation function

$$
 R_k(t_1,t_2) = R_k[f,K](t_1,t_2) = \int_T K(t_1,s,x_{k-1}(s)) \  K(t_2,s,x_{k-1}(s)) \mu(ds) -
$$

$$
\int_T K(t_1,s,x_{k-1}(s)) \ \mu(ds) \cdot \int_T K(t_2,s,x_{k-1}(s)) \ \mu(ds).
$$

\vspace{3mm}

 \ The relation (\ref{final representation} ) in the case $ \ k = m \ $ implies exactly the CLT in the space $ \ C(T,d) \ $ for the
 Monte - Carlo approximation $ \ x_m^m(t) \ $ for the $ \ m^{th} \  $  iteration. \par

\vspace{4mm}

 \ It remains to ground the {\it applicability}  of the CLT for our integrals $ \ x_k^k(t), \ k = 1,2,\ldots,m;  \ t \in T, \ $
 especially for the extremal case $ \ k = m. \ $  We start as  above from the case $ \ x_1(t).\ $ \par

  \ Introduce the following semi - distance on the set $ \ T \ $ depending on some numerical parameter $ \ p(1), \ p(1) \ge 2 \ $

 \begin{equation}  \label{first dist}
 d_{p(1)}(t_1,t_2) := \left\{ \int_T \left[ |K(t_1,s,f(s)) - K(t_2,s,f(s))|^{p(1)} \right] \ \mu(ds)  \right\}^{1/p}.
 \end{equation}

 \ As we know, if

$$
\sup_{t \in T} \int_T |K(t,s,f(s))|^{p(1)} \ \mu(ds) < \infty,
$$

$$
\sup_{t \in T} R_1(t,t) < \infty,
$$
and

$$
\int_0^1 N^{1/p(1)} \left( \ T,d_{p(1)}(\epsilon) \ \right) \ \mu(d \epsilon) < \infty,
$$
then the CLT  for $ \ x_1^1(t) \ $ in the space $ \ C(T,d_{p(1)} \ $  holds true. \par

\vspace{3mm}

 \ Let us consider a general case  $ \ k = 2,3,\ldots,m \ $ by means of induction.
 Introduce as before  the following semi - distances on the set $ \ T \ $ depending on some numerical parameter $ \ p(k), \ p(k) \ge 2 \ $

 \begin{equation}  \label{k th dist}
 d_{p(k)}(t_1,t_2) := \left\{ \int_T \left[ |K(t_1,s,x_{k-1}(s)) - K(t_2,s,x_{k-1}(s))|^{p(k)} \right] \ \mu(ds)  \right\}^{1/p(k)}.
 \end{equation}

 \ We deduce  as above that if

\begin{equation} \label{cond k}
\sup_{t \in T} \int_T |K(t,s,x_{k-1}(s))|^{p(k)} \ \mu(ds) < \infty,
\end{equation}

\begin{equation} \label{bcov}
\sup_{t \in T} R_k(t,t) < \infty,
\end{equation}

and

\begin{equation} \label{k th entr int}
\int_0^1 N^{1/p(k)} \left( \ T,d_{p(k)}, \epsilon \ \right) \ d \epsilon < \infty,
\end{equation}
then the CLT  for $ \ x_k^k(t) \ $ in the space $ \ C(T,d_{p(k)}) \ $  holds true. \par

\vspace{4mm}

\ To summarize: \par

\vspace{4mm}

\ {\bf Theorem 3.1.} Suppose that all the conditions (\ref{cond k}), (\ref{k th entr int}) and (\ref{k th entr int}) are satisfied  for all the
values $ \ k = 1,2,\ldots,m; \ $ as well as all the restrictions formulated before. Assume also that all
the introduced before  distance functions $ \ d_{p(k)}(t_1,t_2)  \ $  are continuous relative the source distance $ \ d = d(t_1,t_2). \ $ \par
\  Then the sequence of random fields $ \ x_m^m[N](t) \ $ satisfies the CLT as $ \ N \to \infty \ $ in the space $ \ C(T,d): \ $

\begin{equation} \label{Theorem 3.1.}
\Law \left\{\sqrt{q[N](m)} \left( \ x_m^{m[N]}(t) - x_m(t) \right) \ \right \} \stackrel{distr}{\to} N(0, R_m(t,s)).
\end{equation}

\vspace{5mm}

\section{The case of non - linear Volterra's integral equations.}

\vspace{5mm}

 \ We consider in this section the non - linear {\it  Volterra's } integral equation  of the ordinary form

 \begin{equation} \label{Volterra eq}
 X(\tau, y) = f(\tau, y) + \int_0^{\tau} \ d \nu \int_T K(\tau, y,\nu, v, X(\nu,v)) \ \mu(dv), \ \tau,\nu \in [0,1],
 \end{equation}
or equally

\begin{equation} \label{Volterra eq const bound}
 X(\tau, y) = f(\tau,y) + \tau \ \int_0^1 d \nu \int_T K(\tau,y, \tau \nu, v, X(\tau \nu,v)) \ \mu(dv).
 \end{equation}

 \ For brevity,

\begin{equation} \label{Volterra  brief}
 X(\tau,y) = f(\tau,y) + F[X](\tau, y),
\end{equation}
where again   $ \ \tau, \nu \in [0,1],  \ y,v \in T, \ $

\begin{equation} \label{FX}
 F[X](\tau, y)= \int_0^{\tau} d\nu \ \int_T K(\tau,y,\nu,v, X(\nu,v)) \ \mu(dv)=
\end{equation}

\begin{equation} \label{FX change}
 \tau \int_0^{\tau} d \nu \int_T K(\tau,y, \tau \nu,v, X(\tau \nu,v)) \ \mu(dv).
\end{equation}

\vspace{4mm}

 \ We retain all the previous notations and restrictions: $ \ (T,d), \ \mu, \ \{Q\} \ $ etc.\par

\vspace{4mm}

 \ The particular case of the equation

 $$
 X(\tau) = X_0 + \int_0^{\tau}  K(\nu,X(\nu)) d\nu
 $$
correspondent to the well - known Cauchy problem for an ordinary differential equation

$$
dX(\tau)/d \tau = K(\tau, X(\tau)), \ \hspace{4mm} X(0) = X_0 = \const.
$$

\vspace{4mm}

\ Let us introduce also the sequence  $ \ \{\eta_i\}, \ i = 1,2,\ldots,N; \ \eta:= \eta_1 $ of  independent uniformly distributed
on the unit interval $ \ [0,1] \ $ r.v., independent also on the source sequence $ \ \{\xi_j\}: $

$$
{\bf P}(\eta_i < w) =  {\bf P}(\eta < w) =  w, \ w \in [0,1],
$$
and of course

$$
{\bf P}(\eta_i < w) =  0, \ w < 0;  \hspace{4mm} {\bf P}(\eta_i < w) =1, \ w > 1.
$$

 \ Note that

\begin{equation} \label{Volterr repres}
F[X](\tau, y) = \tau  \ {\bf E} K(\tau, \ y, \eta \ \tau, \xi, \  X(\tau \ \eta, \ \xi))
\end{equation}
with correspondent consistent  as $ \ n \to \infty \ $ with probability one in the uniform norm $ \  C([0,1] \otimes T)  \ $
Monte - Carlo estimation

\begin{equation} \label{MC approximation}
\hat{F}_n(\tau,y) \stackrel{def}{=} n^{-1} \sum_{i=1}^n \tau \ K(\tau, \ y, \ \eta_i \ \tau, \xi_i, \ X(\tau \eta_i, \ \xi_i)).
\end{equation}

 \vspace{5mm}

  \ We suppose as above that the function $ \ f = f(\cdot,\cdot) \ $ and the
  kernel $ \ K(\cdot,\cdot,\cdot,\cdot,\cdot) \ $ are continuous, bounded and that the kernel $ \ K \ $ satisfies the
 famous Lipschitz condition

\begin{equation}  \label{Lipschitz}
\exists C \in (0,\infty) \ \Rightarrow  | K(\tau,y,\nu,v,z_1) -  K(\tau, y,\nu,v,z_2)| \le C \ |z_1 - z_2|,
\end{equation}
{\it  but we do not assume in contradiction to the previous sections that} $ \ C < 1. \ $\par

\vspace{4mm}

 \ The applications of these equations are described  in particular in the works \cite{Bazarbekov}, \cite{Grigorjeva  Ostrovsky},
 \cite{Golberg}, \cite{Krasnoselskii}, \cite{Tricomi}. \par
  \ The case of solution of super - linear growth  for the kernel $ \ K =  K(\tau,y,\nu,v,z_1) \ $ having in particular blow - up
  solution  is considered in   \cite{Appleby}; see also the reference therein. \par

 \vspace{4mm}

 \ It is no hard to derive by means of induction the following estimation

\begin{equation} \label{good estim}
|F^n[X_1](\tau) - F^n[X_2](\tau)| \le \frac{C^n \ \tau^n}{n!} \ ||X_1 - X_2|| \le  \frac{C^n}{n!} \ ||X_1 - X_2||, \ n = 1,2,\ldots,
\end{equation}
where $ \ F^n \ $ denotes the $ \ n^{th} \ $ iteration of the operator $ \ F, \ $
therefore the {\it continuous} solution of the equation (\ref{Volterra eq}) there exists, is unique and may be
calculated by means of iterations: $ \ X_0 = f(\tau,y), \ $

\begin{equation} \label{Volt iter}
X_{n+1}(\tau,y) = f(\tau,y) + F[X_n](\tau,y),
\end{equation}
see e.g.  \cite{Bazarbekov}, \cite{Golberg}, \cite{Jameson}, \cite{Tricomi}. \par

 \ Error estimate:

\begin{equation} \label{error Volterra}
||X_m(\tau) - X(\tau)||C(T)  \le ||X_1(\tau) - X_0(\tau)||C(T) \cdot \sum_{n=m}^{\infty} \frac{C^n \ \tau^n}{n!} \to 0,
\end{equation}
when $ \ m \to \infty. \ $ \par

\ As a consequence:

\begin{equation} \label{error Volterra sup t}
\sup_{\tau \in [0,1]} \ ||X_m(\tau) - X(\tau)||C(T)  \le  \sup_{\tau \in [0,1] } ||X_1(\tau) - X_0(\tau)||C(T) \cdot \sum_{n=m}^{\infty} \frac{C^n }{n!}.
\end{equation}

\vspace{3mm}

 \ Each integral appearers in (\ref{Volt iter})  may be computed as before by means of depending trial Monte - Carlo method,
so that the rate of convergence in the uniform norm  is equal the classical value $ \ 1/\sqrt{N}. \ $ \par
 \ In detail, we retain the notations and conditions of the second section, especially the vector $ \ \vec{\gamma} \ $ and
 partition of the whole set $ \ S(N), \ $ \hspace{4mm}  $ \ Q = \{Q_k\}, \ k = 1,2,  \ldots,m \ $  (\ref{partition one}),  (\ref{partition general}). \par

 \ Define as above $ \ X_0^0(\tau,y) = X_0(\tau,y) := f(\tau,y) \ $ and recursively

 \begin{equation} \label{first Volterr rand recursion}
X_1^1(\tau,y) = f(\tau,y) +  \frac{\tau}{n(1)} \sum_{i=1}^{n(1)} K(\tau, \ y, \tau \ \eta_i, \ \xi_i, \ X_0^0( \tau \ \eta_i, \ \xi_i)),
 \end{equation}

 \begin{equation} \label{ second rand Volterr recursion}
X_2^2(\tau,y) = f(\tau,y) +  \frac{\tau}{n(2) - n(1)} \sum_{i=n(1)+1}^{n(2)}  K(\tau, \ y, \ \tau \ \eta_i, \ \xi_i, \ X_1^1( \tau \ \eta_i,\xi_i)),
 \end{equation}
and for the values $ \ k = 3,4,\ldots,m \hspace{5mm}  X_k^k(\tau,y) = $

 \begin{equation} \label{ k  rand Volterr  recursion}
f(\tau,y) +  \frac{\tau}{n(k) - n(k-1)} \sum_{i=n(k-1)+1}^{n(k)}  K(\tau, \ y,\ \tau \ \eta_i, \ \xi_i, \ X_{k-1}^{k-1}( \tau \ \eta_i,\xi_i)),
 \end{equation}
and ultimately  $ \ X_m^m(\tau,y) = $

 \begin{equation} \label{m rand Volterr recursion}
X_m^m(\tau,y) = f(\tau,y) +  \frac{\tau}{n(m) - n(m-1)} \sum_{i=n(m-1)+1}^{n(m)}  K(\tau, \ y, \ \tau \ \eta_i, \ \xi_i, \ X_{m-1}^{m-1}(\tau \ \eta_i,\xi_i)).
 \end{equation}

\vspace{5mm}

 \ Recall  once again that $ \ n(m) = N. \ $ \par

\vspace{5mm}

  \ {\sc  Again  preliminary considerations.} \par

\vspace{4mm}

 \ Let us investigate the  first random approximation (\ref{first Volterr rand recursion}). Assume that the r.f
 $ \  (\tau,v) \to \tau \ K(\tau,y, \tau \eta_i, \xi_i, f(\eta_i, \xi_i)) \ $
 satisfies the CLT  in the space $ \ C(([0,1] \otimes T),d); \ $ then one can write (approximately)

$$
X_1^1(\tau,y) = X_1(\tau,y) + \frac{\tau}{n(1)} \sum_{i= 1}^{n(1)} K(\tau, y, \tau \eta_i, \xi_i, X_0^0(\tau \eta_i, \xi_i)) =
$$

 \begin{equation} \label{Volterra recursion 1}
 X_1(\tau,y) +  \frac{1}{\sqrt{ n(1)}} \beta_1(\tau,y) = X_1(\tau,y) + \frac{1}{\sqrt{ q[N](1) }} \beta_1(\tau,y),
 \end{equation}
where $ \ \beta_1(\tau,y) \ $ is a random field such that

\begin{equation} \label{first Volterra asympt1}
\Law(\beta_1 (\cdot,\cdot)) \stackrel{dist}{\to} N \left(0, \ R_1^V[K,f](\tau_1, \tau_2, y_1, y_2) \ \right),
\end{equation}
and

\begin{equation} \label{cov fun V1}
R_1^V[K,f]( \tau_1, \tau_2,y_1,y_2) = R_1^V(\tau_1,\tau_2,y_1,y_2) :=
\end{equation}

\begin{equation} \label{cov fun V2}
\tau_1 \ \tau_2 \  \int_0^1 d \nu  \int_T \ K(\tau_1, y_1, \tau_1 \nu,v, X_0^0(\tau_1 \nu,v)) \ K(\tau_2, y_2, \tau_2 \nu,v, X_0^0(\tau_2 \nu,v)) \mu(dv) -
\end{equation}

\begin{equation} \label{cov fun V3}
\tau_1 \ \int_0^1 d\nu \int_T \ K(\tau_1, y_1, \tau_1 \nu,v, X_0^0(\tau_1 \nu,v)) \mu(dv) \times
\end{equation}

\begin{equation} \label{cov fun V4}
\tau_2 \  \int_0^1 d\nu  \int_T \ K(\tau_2, y_2, \tau_2 \nu, v, X_0^0(\tau_2 \nu,v)) \ \mu(dv).
\end{equation}

\vspace{4mm}

 \ Let us suppose temporarily as above

\begin{equation} \label{first  normal Volt asympt1}
\Law(\beta_1(\cdot)) = N \left(0, R_1^V[f](\tau_1, \tau_2, y_1,y_2) \ \right).
\end{equation}

\vspace{4mm}

 \ The case of the second iteration is more complicated. Let us impose in addition  the following
condition  on the kernel function $ \ K = K(t,s,y,v,z): $

$$
\ \exists \Lambda = \Lambda(\tau, \nu, y,v), \ \exists  \Theta = \Theta(\tau, \nu; y,v, z_1,z_2), \  \exists \delta = \const  \in (0,1]
$$
such that

\begin{equation} \label{restiction}
\sup_{\tau, \nu \in [0,1]} \ \sup_{y,v \in T} |\Lambda(\tau,\nu, y,v)| < \infty,  \ \sup_{\tau, \nu \in [0,1]} \sup_{ z_1,z_2 \in R } |\Theta(\tau,\nu; y,v, z_1,z_2)| < \infty,
\end{equation}

and

$$
 K(\tau, \nu, y,v,z_1) - K(\tau,\nu,y,v,z_2)  =
$$

\begin{equation} \label{smooth  kernel Volt cond}
 \Lambda(\tau, \nu,y,v) \ (z_1-z_2) + \Theta(\tau, \nu ;y, v,z_1,z_2) |z_1-z_2|^{1 + \delta}.
\end{equation}

 \ We derive under this condition (\ref{smooth  kernel Volt cond}) as $ \ N \to \infty \ $ as before

$$
X_2^2(\tau,y) = f(\tau,y) + \frac{1}{n(2) - n(1)} \sum_{i=n(1) + 1}^{n(2)} K(\tau, \tau \eta_i, y, \xi_i, X_1^1(\eta_i \tau, \xi_i))  \sim  f(\tau,y) +
$$

$$
\frac{1}{\gamma(2) \ N } \sum _{i=n(1) + 1}^{n(2)} K\left[(\tau, \tau \eta_i, y, \xi_i, X_1( \tau \eta_i, \xi_i) + [n(1)]^{-1/2} \ \beta_1(\tau \eta_i,\xi_i) )\ \right] \sim
$$

$$
f(\tau,y) +  \tau \int_0^1 d \nu \ \int_T K(\tau, \nu, y,v,  X_1(\tau \nu,v)) \ \mu(d \nu) + \frac{1}{\sqrt{n(2) - n(1)}} \beta_2(\tau,y) =
$$

$$
X_2(\tau,y) + ( q[N](2))^{-1/2}\beta_2(\tau,y),
$$
where $ \ \beta_2(\tau,y) = \beta_2[K,f;N](\tau,y) \ $ is (approximately) Gaussian centered continuous random field

\begin{equation} \label{first Volterra asympt2}
\Law(\beta_2 (\cdot,\cdot)) \stackrel{dist}{\to} N \left(0, \ R_2^V[K,f](\tau_1,\tau_2,y_1,y_2) \ \right),
\end{equation}
with covariation function

\begin{equation} \label{cov fun V2}
R_2^V[K,f](\tau_1, \tau_2, y_1,y_2) = R_2^V(\tau_1, \tau_2, y_1, y_2) :=
\end{equation}

\begin{equation} \label{cov fun V22}
\tau_1 \ \tau_2 \  \int_0^1 d \nu \int_T \ K(\tau_1, y_1, \tau_1 \nu, v, X_1(\tau_1 \nu, v)) \ K(\tau_2, y_2, \tau_2 \nu, v, X_1(\tau_2 \nu,v)) \mu(dv) -
\end{equation}

\begin{equation} \label{cov fun V33}
\tau_1 \ \int_0^1 d \nu \int_T \ K(\tau_1, y_1, \tau_1 \nu, v, X_1(\tau_1 \nu,v)) \mu(dv) \times
\end{equation}

\begin{equation} \label{cov fun V44}
\tau_2 \  \int_0^1 d \nu \int_T \ K(\tau_2, y_2, \tau_2 \nu, v, X_1(\tau_2 \nu,v)) \ \mu(dv).
\end{equation}

\vspace{5mm}

 \ We find quite  analogously for the values $ \ k = 3,4,\ldots,m, \ $ especially for the "ultimate" value $ \ k = m \ $
 the approximate representations $ \ X_k^k(\tau, y) = \ $

\begin{equation} \label{final Vol representation}
 X_k(\tau, y) + (q[N](k))^{-1/2} \beta_k(\tau,y),
\end{equation}
where $ \ \beta_k(\tau, y) = \beta_k[K,f;N](\tau,y) \ $ is as $ \ N \to \infty \ $  approximately Gaussian centered continuous random
field with covariation function

\begin{equation} \label{ k cov fun V2}
R_k^V[K,f](\tau_1, \tau_2, y_1,y_2) = R_k^V(\tau_1, \tau_2, y_1,y_2) := \tau_1 \ \tau_2 \times
\end{equation}

\begin{equation} \label{k cov fun V22}
  \int_0^1 d\nu \int_T \ K(\tau_1, y_1, \tau_1 \nu, v, X_{k-1}(\tau_1 \nu, v)) \ K(\tau_2, y_2, \tau_2 \nu, v, X_{k-1}( \tau_2 \nu, v)) \mu(dv) -
\end{equation}

\begin{equation} \label{k cov fun V33}
\tau_1 \ \int_0^1 d\nu \int_T \ K(\tau_1, y_1, \tau_1 \nu, v, X_{k-1}(\tau_1 \nu, v)) \mu(dv) \times
\end{equation}

\begin{equation} \label{k cov fun V44}
\tau_2 \  \int_0^1 d \nu \int_T \ K(\tau_2, y_2, \tau_2 \nu, v, X_{k-1}(\tau_2 \nu, v)) \ \mu(dv).
\end{equation}

\vspace{5mm}

 \ The relation (\ref{final Vol representation})    in the case $ \ k = m \ $ implies exactly the CLT in the space $ \ C([0.1] \otimes T) \ $ for the
 Monte - Carlo approximation $ \ X_m^m(\tau,y) \ $ for the $ \ m^{th} \  $  iteration. \par

\vspace{4mm}

 \ It remains to ground the {\it applicability}  of the CLT in the space $ \ C([0,1] \otimes T) \ $ for our integrals
 $ \ X_k^k(t,y), \ k = 1,2,\ldots,m;  \ t \in [0,1], \ y \in T \ $
 especially for the extremal case $ \ k = m. \ $  We start as  above from the case $ \ X_1(\tau,y).\ $ \par

  \ Introduce the following semi - distance on the set $ \ U := [0,1] \otimes T \ $ depending on some numerical parameter $ \ p(1), \ p(1) \ge 2 \ $

$$
r_{p(1)}(\tau_1, y_1; \ \tau_2, y_2) \stackrel{def}{=}
$$

 \begin{equation}  \label{first dist}
 \left\{ \int_0^1 d \nu \ \int_T \left[ |K(\tau_1, \nu, y_1, v, f(\nu,v)) - K(\tau_2,\nu, y_2, f(\nu,v))|^{p(1)} \right] \ \mu(d \nu)  \right\}^{1/p(1)}.
 \end{equation}

 \ As we know, if

$$
\sup_{\tau \in [0,1]} \ \sup_{y \in T} \int_0^1 d\nu  \int_T |K(\tau, \nu, y,v, f(\nu,v))|^{p(1)} \ \mu(d\nu) < \infty,
$$

$$
\sup_{\tau \in [0,1]}  \ \sup_{y \in T} R_1^V(\tau,y; \tau,y) < \infty,
$$
and

$$
\int_0^1 N^{1/p(1)} \left( \ U, r_{p(1)}(\epsilon) \ \right) \ d \epsilon < \infty,
$$
then the CLT  for $ \ X_1^1(\tau,y) \ $ in the space $ \ C \left(U,r_{p(1)} \right) \ $  holds true. \par

\vspace{3mm}

 \ Let us consider a general case  $ \ k = 2,3,\ldots,m \ $ by means of induction.
 Introduce as before  the following semi - distances on the set $ \ U \ $ depending on some numerical parameter $ \ p(k) \ $
 greatest or equal than 2: $ \ p(k) \ge 2 \ $

$$
 r_{p(k)}(\tau_1, y_1; \tau_2, y_2) \stackrel{def}{=}
$$

 \begin{equation}  \label{k  th  Volt dist}
 \left\{ \ \int_0^1 d \nu \ \int_T \left[ |K(\tau_1,  \nu, y_1, v, X_{k-1}(\nu,v)) - K(\tau_2,  \nu, y_2, v, X_{k-1}(\nu,v))|^{p(k)} \right] \ \mu(dv) \ \right\}^{1/p(k)}.
 \end{equation}

 \ We deduce  as above that if

\begin{equation} \label{cond Volt k}
\sup_{ \tau \in [0,1]} \sup_{y \in T} \int_0^1 d \nu \ \int_T |K( \tau, \nu, X_{k-1}(\nu,y))|^{p(k)} \ \mu(d\nu) < \infty,
\end{equation}

\begin{equation} \label{bcov}
\sup_{\tau \in [0.1]} \ \sup_{y \in T} R_k^V(\tau, \tau, y,y) < \infty,
\end{equation}

and

\begin{equation} \label{k th entr int}
\int_0^1 N^{1/p(k)} \left( \ U,d_{p(k)}, \epsilon) \ \right) \ d \epsilon < \infty,
\end{equation}
then the CLT  for $ \ X_k^k(t,y) \ $ in the space $ \ C(U, r_{p(k)}) \ $  holds true. \par

\vspace{4mm}

\ To summarize: \par

\ {\bf Theorem 4.1.} Suppose that all the formulated above  conditions  are satisfied  for all the
values $ \ k = 1,2,\ldots,m; \ $ as well as all the restrictions formulated before. Assume also that all
the introduced before  distance functions $ \ r_{p(k)}(t_1,t_2)  \ $  are continuous relative the source distance $ \ d = d(\tau_1, y_1, \tau_2, y_2). \ $ \par
\  Then the sequence of random fields $ \  X_m^m[N](\tau,y) \ $ satisfies the CLT as $ \ N \to \infty \ $ in the space $ \ C(U,d): \ $

\begin{equation} \label{Theorem 4.1}
\Law \left\{\sqrt{q[N](m)} \left( \ X_m^{m[N]}(\tau,y) - X_m(\tau,y) \right) \ \right \} \stackrel{distr}{\to} N(0, R_m^V( \tau_1, \tau_2, y_1,y_2)).
\end{equation}

\vspace{5mm}

\section{ Optimal choice  of partition}.

\vspace{5mm}

 \ Let us  discuss in this section  stated above  the problem of optimal choice  of the partition $ \ Q, \ $  or equally the sequence of
 numbers $ \   \{q(k)\} = \{ q[N](k) \}.  \ $ Recall the restriction (\ref{q norming}):

\begin{equation} \label{q norming again}
\sum_{k=1}^m q[N](k) = n(m) = N
\end{equation}
and a (strong) condition 

\begin{equation} \label{comdition q m strong }
\lim_{N \to \infty} q[N](k) = \infty, \  \lim_{N \to \infty} q[N](m)/N  = 1.
\end{equation}

 \ It follows from the grounding of Theorem 3.3 that under formulated in this theorem conditions $ \ Z = Z_N(q) \stackrel{def}{=}  \Var \left( x_m^{m[N]}(t) \right) \asymp \ $

 \begin{equation} \label{target function}
[q(m)]^{-1} + [q(m) \ q(m-1)]^{-1} + [q(m) \ q(m-1) \ q(m-2)]^{-1} + \ldots +
\end{equation}

\begin{equation} \label{target fun 2}
 [\prod_{k=1}^{m - 1} q(m-k)   ]^{-1}.
 \end{equation}

 \  One can solve  the optimization problem $ \ Z_N(q) \to \min \ $ subject to the  limitation (\ref{q norming again}), as well as taking into account
the positivity of $ \ \{q\} \ $ and restriction (\ref{comdition q m strong }), as ordinary by means of Lagrange's factors method. We will prefer to offer right away
the asymptotically as $ \ N \to \infty \ $  optimal  choice of an optimal  value $ \ \vec{q} = \vec{q_0} \ $ omitting some cumbersome computations:

\begin{equation} \label{optimal q beg}
q_0(m) = q_0[N](m) = \Ent \left[ \ N^{1/2} - C_m N^{1/4}  \ \right],
\end{equation}

\begin{equation} \label{optimal q middle}
q_0(m-1) = q_0[N](m-1) = \Ent \left[ \ N^{1/4} - C_{m-1} N^{1/8} \ \right], \ldots,
\end{equation}

\begin{equation} \label{optimal q end}
q_0(m-k) = q_0[N](m-k) = \Ent \left[ \ N^{ 2^{-k -1}} - C_{m-k} N^{  2^{-k - 2}} \ \right], \
\end{equation}
for the values $ \  k = 2,3,\ldots,m-1; \ $
end finally $ \ q_0(1) = q_0[N](1) = \Ent \left[ \  C_1 N^{2^{-m}} \ \right],  \  $   where $ \ \Ent[z] \ $  denotes the integer part of the real number $ \ z \ $
and $ \ \{C_k\} \ $ are appropriate positive constants. \par

 \ Note that

$$
\forall k = m,m-1,\ldots,1 \ \Rightarrow \lim_{N \to \infty} q_0[N](k) = \infty
$$
and especially

$$
\lim_{N \to \infty} \frac{q_0(m)}{\sqrt{N}} = 1.
$$

\vspace{3mm}

 \ Therefore, all the  conditions of Theorem 3.1 as well as of Theorem 4.1 are satisfied and following the rate of convergence of the
 offered methods by using of elaborated in this section partition  $ \ q_0[N] (k) \ $  is equal to the value $ \ 1/\sqrt{N}, \ $ alike in the classical Monte - Carlo method,
 and  one can apply the Central Limit Theorem in the space of all continuous functions for the error estimation. \par

\vspace{5mm}

\section{Concluding remarks.}

\vspace{5mm}

  \hspace{3mm} {\bf Remark A.}  \ The needed for building of confidence region  for $ \ x_m(t) \ $ covariation function $ \ R_m(t_1,t_2)\ $ may be consistently
  estimated alike the function $ \ x_m^m(t) \ $ itself: $ \ R_m(t_1,t_2) \approx \hat{R}_m(t_1,t_2), \ $ where

$$
\hat{R}_m(t_1,t_2) \stackrel{def}{=} N^{-1} \sum_{i=1}^N   K(t_1,\xi_i,x_{m-1}^{m-1}(\xi_i)) \  K(t_2,\xi_i,x_{m-1}^{m-1}(\xi_i)) -
$$

$$
N^{-1}\sum_{i=1}^N K(t_1,\xi_i, x_{m-1}^{m-1}(\xi_i)) \times N^{-1}\sum_{i=1}^N K(t_2,\xi_i, x_{m-1}^{m-1}(\xi_i)).
$$

 \ Analogous approach may be offered for estimation of covariation function for solving by means of the Monte - Carlo method of the Volterra's equation. \par

\vspace{4mm}

 \ {\bf Remark B.}   The case of {\it systems}  of  non - linear integral equations may be investigated quite analogously. \par

\vspace{4mm}

 \ {\bf Remark C.}   It is interest in our opinion to derive the {\it non - asymptotical} confidence domain for $ \ x_m(\cdot) \ $
 in the uniform norm having at the same range $ \ N^{-1/2} \ $ as in the asymptotical case. \par

\vspace{4mm}

\ {\bf Remark D.}  The important for us the distance function $ \ d_{p(k)}(\cdot,\cdot) \ $ in (\ref{k th dist}) may be estimated as follows.
Note first of all that

$$
||x_k|| \le S(f,\rho) \stackrel{def}{=} \frac{||f||}{1 - \rho}.
$$

 Denote by $ \ B(f,\rho) \ $ the ball in the space $ \ C(T) \ $ with the center at origin having radii $ \ S(f,\rho):m \ $

$$
B(f,\rho) = \stackrel{def}{=} \{ \ g, \ g \in B(f,\rho)\cap C(T,d) \  \}.
$$
 Then

 \begin{equation}  \label{k th dist estim}
 d_{p(k)}(t_1,t_2)  \le \sup_{g \in B(f,\rho)} \left\{ \int_T \left[ |K(t_1,s,g(s)) - K(t_2,s,g(s))|^{p(k)} \right] \ \mu(ds)  \right\}^{1/p(k)}.
 \end{equation}

\vspace{4mm}

 \ {\bf Remark E.}   The remainder term appearing in   (\ref{error Volterra}) may be estimated by means of equality

 $$
 E_{1,m}(z) = z^{1 - m} \left( \ e^z - \sum_{k=0}^{m-2} \frac{z^k}{k!}  \ \right), \ m \ge 3,
 $$
where $ \  E_{\alpha,\beta}(z) \ $ is so - called  generalized Mittag - Leffler's  function of the form

$$
E_{\alpha,\beta}(z) = \sum_{k=0}^{\infty} \frac{z^k}{\Gamma(\alpha k + \beta)}, \ \alpha > 0, \beta  > 0,
$$
see  \cite{Ruzhansky}. This function  is entire having an order $ \ 1/\alpha \ $ and type 1. \par

\vspace{6mm}

\vspace{0.5cm} \emph{Acknowledgement.} {\footnotesize The first
author has been partially supported by the Gruppo Nazionale per
l'Analisi Matematica, la Probabilit\`a e le loro Applicazioni
(GNAMPA) of the Istituto Nazionale di Alta Matematica (INdAM) and by
Universit\`a degli Studi di Napoli Parthenope through the project
\lq\lq sostegno alla Ricerca individuale\rq\rq (triennio 2015 - 2017)}.\par

\end{document}